\documentclass[letterpaper, 10 pt, conference]{ieeeconf}  

 \IEEEoverridecommandlockouts                              

\usepackage{mathtools}
\usepackage{epsfig} 
\usepackage{newtxmath} 
\usepackage{amsmath,amsfonts}
\usepackage{mathtools}
\usepackage{hyperref}
\usepackage{setspace}
\usepackage{enumerate}
\usepackage{cite}
\def\BibTeX{{\rm B\kern-.05em{\sc i\kern-.025em b}\kern-.08em
    T\kern-.1667em\lower.7ex\hbox{E}\kern-.125emX}}

\DeclareMathOperator*{\argmax}{arg\,max}
\DeclareMathOperator*{\argmin}{arg\,min}

\newtheorem{assumption}{Assumption}

\newtheorem{theorem}{Theorem}
\newtheorem{lemma}{Lemma}
\newtheorem{corollary}{Corollary}
\title{\LARGE \bf
Data-Driven Bayesian Nonparametric Wasserstein\\ Distributionally Robust Optimization
}

\author{Chao Ning, \IEEEmembership{Member, IEEE} and Xutao Ma 
\thanks{This work was supported in part by the National Natural Science Foundation
of China under Grant 62103264, and in part by Shanghai Pujiang Program under
Grant 21PJ1406200. (\emph{Corresponding author: Chao Ning})}
\thanks{C. Ning and X. Ma are with Department of Automation, Shanghai
Jiao Tong University, Shanghai 200240, China, and Key Laboratory of System
Control and Information Processing, Ministry of Education of China, Shanghai
200240, China (email: chao.ning@sjtu.edu.cn; maxutao2022@sjtu.edu.cn). }
}

\begin{document}

\allowdisplaybreaks

\maketitle
\thispagestyle{empty}
\pagestyle{empty}

\begin{abstract}

In this work, we develop a novel data-driven Bayesian nonparametric Wasserstein distributionally robust optimization (BNWDRO) framework for decision-making under uncertainty. The proposed framework unifies a Bayesian nonparametric method and the Wasserstein metric to decipher the global-local features of uncertainty data and encode these features into a novel data-driven ambiguity set. By establishing the theoretical connection between this data-driven ambiguity set and the conventional Wasserstein ambiguity set, we prove that the proposed framework enjoys the finite sample guarantee and asymptotic consistency. To efficiently solve the resulting distributionally robust optimization problem under the BNWDRO framework, we derive for this optimization problem an equivalent reformulation, which is kept tractable for many practical scenarios. Numerical experiments on a unit commitment problem verify the effectiveness of the proposed BNWDRO framework compared with existing methods.

\end{abstract}


\section{Introduction}

In the context of decision-making under uncertainty, estimating the uncertainty distribution is of the utmost importance. Stochastic programming (SP) is a well-developed method for addressing this issue \cite{shapiro2021lectures}. It either assumes the uncertainty distribution is known or uses an estimated distribution in substitution for the true one. However, the true distribution is not perfectly known in practice, and the distribution estimated from partial information inevitably deviates from the true one. Therefore, the SP approach can lead to poor out-of-sample performance due to the distribution deviation, and this effect is commonly referred to as the "optimizer's curse" \cite{smith2006optimizer}.

In recent years, distributionally robust optimization (DRO) has emerged as a promising data-driven method to hedge against this effect by leveraging ambiguity sets \cite{wiesemann2014distributionally}.
In the literature, ambiguity sets can be broadly classified into moment-based ambiguity sets and discrepancy-based ambiguity sets. To construct a moment-based ambiguity set, the first- and second-order moment information of the data is commonly leveraged \cite{li2021distributionally}. 
References \cite{wiesemann2014distributionally} and \cite{bertsimas2019adaptive} further extended the moment-based ambiguity set to a second-order cone (SOC) representable ambiguity set.
Alternatively, the decision-maker can construct a discrepancy-based ambiguity set by leveraging a discrepancy function to measure the dissimilarity of probability distributions \cite{shapiro2021lectures,10111046,10167717}.
{For example, relative entropy pseudo-distance can be used for constructing a discrepancy-based ambiguity set \mbox{\cite{petersen2000minimax,hu2013kullback,falconi2023distributionally}}. However, such an ambiguity set only contains distributions absolutely continuous with respect to the nominal one, thereby failing to include the true distribution in some cases.
On the contrary, using a probability metric is without such a defect.} In recent years, the Wasserstein metric-based ambiguity set has aroused considerable attention in the DRO research due to its asymptotic convergence and tractable reformulation \cite{mohajerin2018data,gao2023distributionally,chen2022data,doi:10.1137/22M1494105,mcallister2023inherent}, and it has witnessed many applications \cite{zhu2019wasserstein,yang2020wasserstein,lee2018minimax,guo2018data,zhao2023minimax,coulson2021distributionally,10153787}.

In practice, DRO methods with theoretical performance guarantees are favored, such as the Wasserstein DRO method and the moment-based DRO method. However, when confronted with complex data structures, these DRO methods typically overlook the global-local features of uncertainty data, thus leading to conservative solutions. Few research studies have sought to encode such global-local features into the ambiguity set 
\cite{ning2021online,esteban2022partition,hanasusanto2015distributionally}. However, these approaches either lack theoretical performance guarantees or rely on prior knowledge of these features to achieve such guarantees instead of learning from data, turning out to be impractical in real-world scenarios. Therefore, there exists a research gap in developing a data-driven DRO framework that can fully explore and characterize the global-local features of the data, and at the same time enjoy theoretical performance guarantees for decision-making under uncertainty.

To fill the research gap, this work proposes a novel data-driven Bayesian nonparametric Wasserstein DRO (BNWDRO) framework for decision-making under uncertainty. Instead of presuming information on the data structure, we leverage a Bayesian nonparametric method known as the Dirichlet process mixture model (DPMM) to automatically decompose the data into distinct clusters. Based on the results of the DPMM, we construct for each cluster a local Wasserstein ball, thus allowing for an accurate characterization of the local distributional features. Subsequently, a novel data-driven ambiguity set is cast as the weighted Minkovski sum of these local Wasserstein balls. We investigate the connection between the proposed data-driven ambiguity set and the conventional Wasserstein ambiguity set. By this connection, we prove that the proposed BNWDRO framework enjoys both finite sample guarantee and asymptotic consistency. To solve the resulting DRO problem with respect to the proposed data-driven ambiguity set, we derive equivalent reformulations for this optimization problem.
Compared with other DRO methods, the proposed BNWDRO framework fully explores and characterizes the global-local features of uncertainty data, and at the same time, its performance is theoretically guaranteed, thereby yielding high-quality decisions for the decision-maker. 

\textit{Contributions:} (i) We propose a novel BNWDRO framework that deciphers and encodes the global-local features of data into a data-driven ambiguity set. (ii) We establish a theoretical connection between the proposed data-driven ambiguity set and the Wasserstein ambiguity set, based on which we prove that the proposed BNWDRO framework enjoys theoretical performance guarantees. (iii) We derive equivalent reformulations for the resulting DRO problem under our framework. (iv) We demonstrate empirically the effectiveness of the proposed BNWDRO framework on a unit commitment problem (UC) compared with existing methods.
 
\textit{Notation:}
The set of positive running indices up to $N$ is denoted by $[N]$. $\delta_{(\cdot)}$ is the Dirac delta distribution. 
$\Gamma^N$ denotes the set containing $N$ data $\widehat{\boldsymbol{w}}_l,l\in [N]$, and $\widehat{\mathbb{P}}^N=\sum_{l \in [N]}\delta_{\widehat{\boldsymbol{w}}_l}/N$ is the corresponding empirical distribution. For any set $\mathscr{A}$, we denote its cardinality by $|\mathscr{A}|$, and its support function is defined as $\sigma (\mathscr{A})(\boldsymbol{z})=\sup_{\boldsymbol{w}\in \mathscr{A}}=\langle \boldsymbol{z},\boldsymbol{w}\rangle$. $\overline{\mathbb{R}}=\mathbb{R}\cup \{ \infty,-\infty \}$ is the extended real number line. {By $\mathcal{P}^N$ we denote the $N$ fold product of a probability measure $\mathcal{P}$.}

\section{Bayesian Nonparametric Wasserstein Distributionally Robust Optimization Framework}
\label{framework}

Consider the following SP problem
\begin{equation}
\label{SP}
   \min_{x\in \mathscr{X}}\Big\{ f(\boldsymbol{x})+\mathbb{E}_{\mathcal{P}}[g(\boldsymbol{x},\boldsymbol{w})]\Big\}
\end{equation}
where the minimization is taken over a feasible region $\mathscr{X}\subseteq \mathbb{R}^n$, uncertainty $\boldsymbol{w}$ is supported on $\varXi\subseteq \mathbb{R}^m$ with distribution $\mathcal{P}$, and $f(\cdot):\mathbb{R}^n\to \overline{\mathbb{R}}$ and $g(\cdot,\cdot):\mathbb{R}^n\times \mathbb{R}^m\to \overline{\mathbb{R}}$ are loss functions.

However, instead of knowing the true distribution $\mathcal{P}$ perfectly, we only have access to a set $\Gamma^N$ containing $N$ historical uncertainty data $\widehat{\boldsymbol{w}}_l,l\in [N]$. Therefore, we aim to construct a data-driven ambiguity set $\mathscr{B}(\Gamma^N)$ to capture the distributional ambiguity by leveraging information in data.

In constructing the ambiguity set, when confronted with complex data structures, previous research fails to provide a theoretically guaranteed method to capture these unrevealed global-local features of data. 

To fill this research gap, we propose a BNWDRO framework. {In this framework, the first step is to learn features from data. To this end, we adopt a Bayesian nonparametric method known as the DPMM to automatically learn the clustering label of data. In the second step, based on the clustering results, we construct a data-driven ambiguity set that accurately characterizes the local distribution of each cluster.}
\vspace{-6pt}
\subsection{Dirichlet Process Mixture Model Clustering}
\vspace{-5pt}

{The DPMM is a well-known Bayesian nonparametric method for clustering and the number of clusters can be automatically learned from data instead of being given as prior knowledge.
The sampling procedure of the DPMM is presented as follows.}
\vspace{-8pt}
\allowdisplaybreaks
    \begin{gather}
        p\sim \text{DP}(h, H)\nonumber\\
        \theta\sim p\label{DP1}\\
        \boldsymbol{w}\sim F({\phi})\nonumber
    \end{gather}
{where $h,H,F(\phi)$ are hyperparameters of the DPMM, $F(\phi)$ is a distribution family with parameter $\phi$, DP represents the Dirichlet process characterized by a scale parameter $h$ and a prior measure $H$
on the space of parameter {${\phi}$}, and $p$ is a sample from the Dirichlet process and it is also a probability
measure on the space of parameter {${\phi}$}}

To facilitate data clustering, we reformulate (\ref{DP1}) as follows to incorporate the clustering labels as latent variables. 
\begin{equation}
    \label{DP2}
    \begin{aligned}
        \boldsymbol{\gamma}&\sim \text{GEM}(h)\  \ \ \ \  \ \ \ \ \phi_k \sim H\\
        z &\sim \text{Mult}(\boldsymbol{\gamma}) \ \ \ \ \ \ \ \  \boldsymbol{w} \sim F({\phi_{z}})
    \end{aligned}
\end{equation}
where $\boldsymbol{\gamma}=(\gamma_k),k\in \mathbb{N}$ is the mixing proportion sampled from 
the stick-breaking distribution {denoted by} GEM and
 $z$ is the clustering label. 

In relation to (\ref{DP1}), it holds that 
\begin{equation}
    \label{disp}
    p=\sum_{k=1}^{+\infty}
\gamma_k\delta_{\phi_k}
\end{equation}
where $\gamma_k$ represents the weight of component $k$ with {local} distribution $F({\phi_k})$,
and
$\delta_{(\cdot)}$ is the Dirac delta distribution.

{Since no prior information is available to us, we follow the common practice in the DPMM research to take $F(\phi)$ to be the Gaussian distribution family $\mathcal{N}(\boldsymbol{\mu},\boldsymbol{\Sigma})$
with its parameter $\phi$=$(\boldsymbol{\mu},\boldsymbol{\Sigma})$.
The prior distribution of $(\boldsymbol{\mu},\boldsymbol{\Sigma})$ is set to be a normal Wishart distribution
NW$(\boldsymbol{\mu}_0,\lambda_0,\boldsymbol{W}_0,\nu_0)$.}
Supposing that we now have data $\widehat{\boldsymbol{w}}_l,l\in [N]$ independently sampled from (\ref{DP2}), the next step is to estimate the 
posterior expected value of $\gamma_k,\phi_k$ in (\ref{disp}) as well as the latent clustering label $z_l$ associated with the
data $\widehat{ \boldsymbol{w}}_l$. 
{To this aim, we adopt the well-developed variational inference approach \mbox{\cite{blei2006variational}}, which features a fast computational speed and has been integrated into many tool packages such as \mbox{\emph{sklearn}}.}

{By the variational inference, the posterior distribution of the latent label $z_l$ is derived, and we further define clustering label $\overline{z_l}$ associated with data $\widehat{\boldsymbol{z}}_l$ as follows}
\begin{equation}
    {\overline{z_l}=\argmax_{k} p(z_l=k|\boldsymbol{w}_1,\cdots,\boldsymbol{w}_N)}
\end{equation}
{According to the clustering label $\overline{z_l}$, the original data set $\Gamma^N$ is partitioned into $K$ clusters $\Gamma^N_k,k\in [K]$.}

\vspace{-3pt}

\subsection{Ambiguity Set Construction}
\vspace{-4pt}
{Based on the clustering results, we construct for each cluster a local Wasserstein ball $\mathscr{W}_{\theta_k} ( \widehat{\mathbb{P}}_k^N )$ centered at the local empirical distribution $ \widehat{\mathbb{P}}_k^N =\sum_{\widehat{\boldsymbol{w}}\in \Gamma^N_k}\delta_{\widehat{\boldsymbol{w}}}/|\Gamma^N_k|$ as follows.}
\allowdisplaybreaks
\begin{gather}
        \mathscr{W}_{\theta_k} \left( \widehat{\mathbb{P}}_k^N \right)=\left\{ \mathbb{Q}\in \mathscr{P}(\varXi)\left| d_W\left( \mathbb{Q},\widehat{\mathbb{P}}_k^N \right)
        \leq \theta_k \right. \right\}\\
       d_W\left( \mathbb{Q},\widehat{\mathbb{P}}_k^N\right)=\inf_{\pi \in \Pi(\mathbb{Q},\widehat{\mathbb{P}}_k^N)}\int_{\varXi \times \varXi}\| \boldsymbol{x}-\boldsymbol{y} \| d\pi(\boldsymbol{x},\boldsymbol{y})
\end{gather}
where $ \mathscr{P}(\varXi)$ is the set of all Borel probability measures
supported on $\varXi$ , $\theta_k$ is the Wasserstein radius, $d_W( \cdot,\cdot )$ is the Wasserstein metric, $\| \cdot \|$ is a norm on $\mathbb{R}^m$,
and $\Pi(\mathbb{Q},{\mathbb{P}})$
contains all measures with marginal distributions $\mathbb{Q}$ and ${\mathbb{P}}$. For more information about the Wasserstein ambiguity set, we refer to \cite{mohajerin2018data}.

Unifying these local Wasserstein balls, we construct a novel data-driven ambiguity set $\mathscr{B}(\Gamma^N)$ as follows.
\begin{equation}
\label{def_amb}
    \mathscr{B}(\Gamma^N)=\oplus_{k=1}^K\ \frac{|\Gamma^N_k |}{| \Gamma^N |} \mathscr{W}_{\theta_k} \left( \widehat{\mathbb{P}}_k^N  \right)
\end{equation}
where $\oplus$ represents the Minkowski sum.

This data-driven ambiguity set can accurately capture the local distributional features of each cluster, and at the same time, the global structure of these clusters is precisely reflected.

Based on this novel data-driven ambiguity set $\mathscr{B}(\Gamma^N)$, we formulate the BNWDRO counterpart of SP problem (\ref{SP}) as follows.
\begin{equation}
\label{BNWDRO}
       \min_{x\in \mathscr{X}} \bigg\{ f(\boldsymbol{x})+\max_{\mathbb{Q}\in \mathscr{B}(\Gamma^N)}\mathbb{E}_{\mathbb{Q}}[g(\boldsymbol{x},\boldsymbol{w})]\bigg\} 
\end{equation}

\section{Theoretical Results}

In this section, we first develop the theoretical performance guarantee of the proposed BNWDRO framework. Subsequently, the equivalent reformulation of the resulting BNWDRO problem (\ref{BNWDRO}) is derived.

\subsection{Finite Sample Guarantee and Asymptotic Consistency}

We begin by establishing the connection between the proposed data-driven ambiguity set $\mathscr{B}(\Gamma^N)$ and the Wasserstein ambiguity set $\mathscr{W}_{\theta} ( \widehat{ \mathbb{P} }^N  )$, given by the following theorem.

\begin{theorem}
\label{th1}
    Let $\underline{\theta}^N=\min_{k\in[K]} \frac{|\Gamma^N_k|}{|\Gamma^N|}\theta_k$ and $\overline{\theta}^N=\max_{k\in[K]} \theta_k$, and then it holds that 
    \begin{equation}
        \mathscr{W}_{\underline{\theta}^N}\Big( \widehat{ \mathbb{P} }^N  \Big) \subseteq \mathscr{B}(\Gamma^N) \subseteq \mathscr{W}_{\overline{\theta}^N}\Big( \widehat{ \mathbb{P} }^N \Big) 
    \end{equation}
\end{theorem}

\vspace{7pt}
\begin{proof}
 {According to the disintegration theorem, $\pi \in \Pi(\mathbb{Q},\widehat{\mathbb{P}}_k^N)$ can be expressed as}
     \begin{equation}
         {\pi=\widehat{\mathbb{P}}^N_k \otimes \mathbb{Q}_{w}=\frac{1}{|\Gamma^N_k|}\sum_{\widehat{w}\in \Gamma^N_k} \delta_{\widehat{w}}\otimes \mathbb{Q}_{{w}}}
     \end{equation}
     {where $\mathbb{Q}_{w}$ is a transition kernel from $\varXi$ to $\varXi$ and the input variable is expressed explicitly as subscription $w$. Therefore, the marginal probability measure $\mathbb{Q}$ is the weighted sum of conditional distribution $\mathbb{Q}_{\widehat{w}}$ as follows.}
     \begin{equation}
        { \mathbb{Q}=\pi(\varXi,\cdot)=\frac{1}{|\Gamma^N_k|}{\sum_{\widehat{w}\in \Gamma^N_k} }\mathbb{Q}_{\widehat{w}}}\label{r2q30}
     \end{equation}

{Based on (\mbox{\ref{r2q30}}), the Wasserstein ambiguity set centered at local empirical distribution $\widehat{\mathbb{P}}_k^N$ with radius $\theta_k$ admits the following reformulation.}
\begin{equation}
\label{local_w}
    \mathscr{W}_{\theta_k} \left( \widehat{\mathbb{P}}_k^N  \right)=\left\{ \mathbb{Q} \left|
    \begin{gathered}
        \mathbb{Q}=\sum_{\widehat{\boldsymbol{w}}\in \Gamma^N_k}\frac{1}{| \Gamma^N_k |} \mathbb{Q}_{\widehat{\boldsymbol{w}}},\\
        \mathbb{Q}_{\widehat{\boldsymbol{w}}} \in \mathscr{P}(\varXi),\forall \ \widehat{\boldsymbol{w}}\in \Gamma^N_k,\\
        \sum_{\widehat{\boldsymbol{w}}\in \Gamma^N_k} \frac{1}{| \Gamma^N_k |} \mathbb{E}_{\boldsymbol{w}\sim \mathbb{Q}_{\widehat{\boldsymbol{w}}}}  \left\| \boldsymbol{w}-\widehat{\boldsymbol{w}} \right\|\leq \theta_k
    \end{gathered}   \right. \right\}
\end{equation}

\begin{spacing}{0}
Based on the reformulation (\ref{local_w}) for the local ambiguity sets, we can reformulate the data-driven ambiguity set $\mathscr{B}(\Gamma^N)$ as follows.
\begin{gather}
    \mathscr{B}(\Gamma^N)=\left\{ \mathbb{Q}\left|\mathbb{Q}=\sum_{i=k}^K \frac{|\Gamma^N_k |}{| \Gamma^N |} \mathbb{Q}_k, \mathbb{Q}_k \in\mathscr{W}_{\theta_k} \left( \widehat{\mathbb{P}}_k^N  \right)   \right. \right\}\\
    =\left\{   \mathbb{Q}\left|
    \begin{gathered}
    \mathbb{Q}=\sum_{i=k}^K \frac{|\Gamma^N_k |}{| \Gamma^N |} \left[  \sum_{\widehat{\boldsymbol{w}}\in \Gamma^N_k}\frac{1}{| \Gamma^N_k |} \mathbb{Q}_{\widehat{\boldsymbol{w}}} \right], \\
    \mathbb{Q}_{\widehat{\boldsymbol{w}}}\in \mathscr{P}(\varXi),\forall\ \widehat{\boldsymbol{w}}\in \Gamma^N, \\
            \sum_{\widehat{\boldsymbol{w}}\in \Gamma^N_k} \frac{1}{| \Gamma^N_k |} \mathbb{E}_{\boldsymbol{w}\sim \mathbb{Q}_{\widehat{\boldsymbol{w}}}}  \left\| \boldsymbol{w}-\widehat{\boldsymbol{w}} \right\|\leq \theta_k,\forall k\in [K]
    \end{gathered}
    \right. \right\}
    \end{gather}
    
    \begin{gather}
        =\left\{   \mathbb{Q}\left|
    \begin{gathered}
    \mathbb{Q}= \frac{1}{| \Gamma^N |}   \sum_{\widehat{\boldsymbol{w}}\in \Gamma^N}\mathbb{Q}_{\widehat{\boldsymbol{w}}},  \\
    \mathbb{Q}_{\widehat{\boldsymbol{w}}}\in \mathscr{P}(\varXi),\forall\ \widehat{\boldsymbol{w}}\in \Gamma^N, \\
            \sum_{\widehat{\boldsymbol{w}}\in \Gamma^N_k} \frac{1}{| \Gamma^N|} \mathbb{E}_{\boldsymbol{w}\sim \mathbb{Q}_{\widehat{\boldsymbol{w}}}}  \left\| \boldsymbol{w}-\widehat{\boldsymbol{w}} \right\|\leq \frac{|\Gamma^N_k|}{|\Gamma^N|} \theta_k,\forall k\in [K]
    \end{gathered}
    \right. \right\} \label{final}
\end{gather}
\end{spacing}

\vspace{1.5em}
Note that similar to the local Wasserstein ambiguity set, the Wasserstein ambiguity set centered at the overall empirical distribution $\widehat{\mathbb{P}}^N$ with radius $\theta$ can be reformulated as
\begin{equation}
\label{global_w}
    \mathscr{W}_{\theta} \left( \widehat{ \mathbb{P} }^N  \right)=\left\{ \mathbb{Q} \left|
    \begin{gathered}
        \mathbb{Q}=\sum_{\widehat{\boldsymbol{w}}\in \Gamma^N}\frac{1}{| \Gamma^N |} \mathbb{Q}_{\widehat{\boldsymbol{w}}},\\
        \mathbb{Q}_{\widehat{\boldsymbol{w}}} \in \mathscr{P}(\varXi),\forall \ \widehat{\boldsymbol{w}}\in \Gamma^N,\\
        \sum_{\widehat{\boldsymbol{w}}\in \Gamma^N} \frac{1}{| \Gamma^N |} \mathbb{E}_{\boldsymbol{w}\sim \mathbb{Q}_{\widehat{\boldsymbol{w}}}}  \left\| \boldsymbol{w}-\widehat{\boldsymbol{w}} \right\|\leq \theta
    \end{gathered}   \right. \right\}
\end{equation}

By comparing (\ref{final}) and (\ref{global_w}), we notice that the difference between the Wasserstein ambiguity set $\mathscr{W}_{\theta} ( \widehat{ \mathbb{P} }^N  )$ and the proposed data-driven ambiguity set $\mathscr{B}(\Gamma^N)$ lies in the constraints of the conditional distribution $\mathbb{Q}_{\widehat{\boldsymbol{w}}}$. In the Wasserstein ambiguity set, all the conditional distributions are constrained together, while the proposed data-driven ambiguity set constraints conditional distributions of different clusters separately.
Focusing on this difference, we note that 
\begin{gather}
        \sum_{\widehat{\boldsymbol{w}}\in \Gamma^N} \frac{1}{| \Gamma^N |} \mathbb{E}_{\boldsymbol{w}\sim \mathbb{Q}_{\widehat{\boldsymbol{w}}}}  \left\| \boldsymbol{w}-\widehat{\boldsymbol{w}} \right\|\leq\underline{\theta}^N \\
        \Rightarrow \sum_{\widehat{\boldsymbol{w}}\in \Gamma^N_k} \frac{1}{| \Gamma^N |} \mathbb{E}_{\boldsymbol{w}\sim \mathbb{Q}_{\widehat{\boldsymbol{w}}}}  \left\| \boldsymbol{w}-\widehat{\boldsymbol{w}} \right\|\leq \min_{k\in[K]} \left\{\frac{|\Gamma^N_k|}{|\Gamma^N|}\theta_k \right\}\\
        \Rightarrow \sum_{\widehat{\boldsymbol{w}}\in \Gamma^N_k} \frac{1}{| \Gamma^N |} \mathbb{E}_{\boldsymbol{w}\sim \mathbb{Q}_{\widehat{\boldsymbol{w}}}}  \left\| \boldsymbol{w}-\widehat{\boldsymbol{w}} \right\|\leq \frac{|\Gamma^N_k|}{|\Gamma^N|}\theta_k 
\end{gather}
Therefore, the first set containment relationship holds.

Also, notice that 
\begin{gather}
    \sum_{\widehat{\boldsymbol{w}}\in \Gamma^N} \frac{1}{| \Gamma^N |} \mathbb{E}_{\boldsymbol{w}\sim \mathbb{Q}_{\widehat{\boldsymbol{w}}}}  \left\| \boldsymbol{w}-\widehat{\boldsymbol{w}} \right\|=
 \sum_{k=1}^K \sum_{\widehat{\boldsymbol{w}}\in \Gamma^N_k} \frac{1}{| \Gamma^N |} \mathbb{E}_{\boldsymbol{w}\sim \mathbb{Q}_{\widehat{\boldsymbol{w}}}}  \left\| \boldsymbol{w}-\widehat{\boldsymbol{w}} \right\|\nonumber\\
 \leq \sum_{k=1}^K\frac{|\Gamma^N_k|}{|\Gamma^N|} \theta_k
 \leq \sum_{k=1}^K\frac{|\Gamma^N_k|}{|\Gamma^N|}\max_{k\in[K]}\theta_k=\overline{\theta}^N\label{r2q9}
\end{gather}
Therefore, the second set containment relationship holds.
\end{proof}

With Theorem \ref{th1} bridging the proposed data-driven ambiguity set and the conventional Wasserstein ambiguity set, we can prove the finite sample guarantee and the asymptotic consistency of the proposed BNWDRO framework.

To state the finite sample guarantee and the asymptotic consistency, we define the optimal solution $\widehat{\boldsymbol{x}}_N$ and the optimal value $\widehat{J}_N$ of the BNWDRO problem (\ref{BNWDRO}) below.
        \begin{gather}
         \widehat{\boldsymbol{x}}_N=\argmin_{x\in \mathscr{X}} \{  f( \boldsymbol{x})+ \max_{\mathbb{Q}\in \mathscr{B}(\Gamma^N)}\mathbb{E}_{\mathbb{Q}}[g( \boldsymbol{x},\boldsymbol{w})]\}\label{x1}\\
         {\scriptstyle \widehat{J}_N}= \min_{x\in \mathscr{X}}\{  f( \boldsymbol{x})+ \max_{\mathbb{Q}\in \mathscr{B}(\Gamma^N)}\mathbb{E}_{\mathbb{Q}}[g( \boldsymbol{x},\boldsymbol{w})]\}\label{J1}
    \end{gather}

The following measure concentration in Lemma \ref{lemma1} will be the cornerstone in the proof of the theoretical guarantee.
\begin{assumption}
\label{assum}
    The true distribution $\mathcal{P}$ is a light-tailed distribution, \emph{i.e.}, there existing a constant $b>1$ such that 
    \begin{equation}
         \mathbb{E}_{\mathcal{P}}\left[ {exp}\left( \| \boldsymbol{w} \|^b \right) \right]< \infty
    \end{equation}
\end{assumption}
\vspace{5pt}

\begin{lemma} (Theorem 3.4 \cite{mohajerin2018data})
\label{lemma1}
\setlength{\arraycolsep}{0.5pt}
    Suppose that Assumption \ref{assum} holds. Then the following measure concentration inequality holds for a prescribed $0< \beta <1$.
    \begin{equation}
    \label{r2q41}
         \mathcal{P}^{N}\left\{  d_W( \mathcal{P},\widehat{\mathbb{P}}^N)\geq \epsilon^N(\beta) \right\}\leq \beta
    \end{equation}
    where $\epsilon^N(\beta)$ is determined as follows by $\beta$ and positive constants $b_1,b_2$ which depend on the shape of  $\mathcal{P}$.
    \begin{equation}
         \epsilon^N(\beta)=\left\{  \begin{array}{ll}
            \left( \frac{log(b_1\beta^{-1})}{b_2N} \right)^{1/max\{ m,2 \}}&,\text{if} \ N\geq \frac{log(b_1\beta^{-1})}{b_2} \\
            \left( \frac{log(b_1\beta^{-1})}{b_2N} \right)^{1/a}&,\text{if} \ N\geq \frac{log(b_1\beta^{-1})}{b_2}
        \end{array}  \right.
    \end{equation}
\end{lemma}
\vspace{5pt}

Based on Theorem \ref{th1} and Lemma \ref{lemma1}, we prove the finite sample guarantee and asymptotic consistency of the proposed BNWDRO framework in Theorem \ref{finite} and Theorem \ref{asym}, respectively.

\begin{theorem}
\label{finite}
    (Finite sample guarantee) Suppose that Assumption \ref{assum} holds and that $\underline{\theta}^N\geq\epsilon^N(\beta) $. Then we have the following finite sample guarantee.
    \begin{equation}
    \label{gua}
        \mathcal{P}^{N}\left\{  \widehat{\varXi}_N: f( \widehat{\boldsymbol{x}}_N)+ \mathbb{E}_{\mathcal{P}}[g( \widehat{\boldsymbol{x}}_N,\boldsymbol{w})]\leq  \widehat{J}_N \right\}\geq 1-\beta
    \end{equation}
    where $\beta\in (0,1)$ is a prescribed significance parameter and $\widehat{J}_N$ serves as a certificate of the out-of-sample performance of decision $\widehat{\boldsymbol{x}}_N$.
\end{theorem}
\vspace{5pt}
\begin{proof}
    (\ref{gua}) holds because
    \begin{equation*}
    \begin{gathered}
        \mathcal{P}^{N}\left\{  \widehat{\varXi}_N: f( \widehat{\boldsymbol{x}}_N)+ \mathbb{E}_{\mathcal{P}}[g( \widehat{\boldsymbol{x}}_N,\boldsymbol{w})]\leq  \widehat{J}_N \right\}
        \geq \mathcal{P}^{N}\left\{ \mathcal{P}\in \mathscr{B} \right\}  \\\geq \mathcal{P}^{N}\left\{ \mathcal{P}\in \mathscr{W}_{\underline{\theta}^N}(\widehat{\mathbb{P}}^N) \right\}
        \geq 1-\mathcal{P}^{N}\left\{  d_W( \mathcal{P},\widehat{\mathbb{P}}^N)\geq \underline{\theta}^N \right\}\\
        \geq1-\mathcal{P}^{N}\left\{  d_W( \mathcal{P},\widehat{\mathbb{P}}^N)\geq \epsilon^N(\beta)  \right\}\geq 1-\beta
    \end{gathered}
    \end{equation*}
    where the second inequality comes from Theorem \ref{th1} and the last inequality comes from Lemma \ref{lemma1}.
\end{proof}

\begin{theorem}
\label{asym}
    (Asymptotic consistency) Suppose that Assumption \ref{assum} holds and that $\overline{\theta}^N\leq \epsilon^{N}(\beta^N)$ for a sequence of $\beta^N\in (0,1),N\in \mathbb{N}$, which satisfies $\lim_{N\to \infty}\epsilon^{N}(\beta^N)=0$ and $\sum_{N=1}^{\infty}\beta^N\leq +\infty$. Let $\widehat{\boldsymbol{x}}_N$ and $\widehat{J}_N$ be defined in (\ref{x1}) and (\ref{J1}), respcetively. Then, the following asymptotic consistency holds.
    \begin{enumerate}[i)]
        \item If $g(\boldsymbol{x},\boldsymbol{w})$ is continuous in $\boldsymbol{w}$ for all $\boldsymbol{x}\in \mathscr{X}$ and $|g(\boldsymbol{x},\boldsymbol{w})|$ is upper bounded by $L(1+\| \boldsymbol{w} \|)$ with respect to a constant $L\geq 0$ for all $\boldsymbol{x}\in \mathscr{X}$ and $\boldsymbol{w}\in \varXi$, then $\mathcal{P}^{\infty}$-almost surely $\lim_{N\to \infty}\widehat{J}_N=J^\ast$, where $J^\ast$ is the optimal value of the SP problem (\ref{SP}).
        \item If the assumption in i) hold, the region $\mathscr{X}$ is a closed set, and for each $\boldsymbol{w}\in \varXi$, $g(\boldsymbol{x},\boldsymbol{w})$ is continuous in $x$, then any accumulation point of $\{ \widehat{\boldsymbol{x}}^N :N\in \mathbb{N}\}$ is an optimal solution to the SP problem (\ref{SP}) $\mathcal{P}^{\infty}$-almost surely.
    \end{enumerate}
\end{theorem}

\begin{proof}
    According to Lemma 3.7 in \cite{mohajerin2018data}, for any sequence $\mathbb{Q}^N\in \mathscr{W}_{\epsilon^{N}(\beta^N)}(\widehat{\mathbb{P}}^N)$, we have 
    \begin{equation}
        d_W(\mathbb{Q}^N,\mathcal{P})\xrightarrow{{{N\to \infty}}}0,\ \mathcal{P}^{\infty}-\text{a.s.}
    \end{equation}
    
    {By Theorem \mbox{\ref{th1}}, $\mathscr{B}(\Gamma^N)$ is a subset of $\mathscr{W}_{\epsilon^{N}(\beta^N)}(\widehat{\mathbb{P}}^N)$. Therefore, for any sequence $\tilde{\mathbb{Q}}^N\in \mathscr{B}(\Gamma^N)$ we also have}
        \begin{equation}
        \label{asw}
        d_W(\tilde{\mathbb{Q}}^N,\mathcal{P})\xrightarrow{{N\to \infty}}0,\ \mathcal{P}^{\infty}-\text{a.s.}
    \end{equation}
    
By leveraging (\ref{asw}), we can prove i) and ii) by the same argument as in the proof of Theorem 3.6 in \cite{mohajerin2018data}.
\end{proof}

\vspace{-3pt}
\subsection{Tractable Reformulation}
\vspace{-4pt}
Notice that the worst-case expectation in (\ref{BNWDRO}) can be reformulated as (\ref{mid1}).
\begin{gather}
    \max_{\mathbb{Q}\in \mathscr{B}(\Gamma^N)}\mathbb{E}_{\mathbb{Q}}[g(\boldsymbol{x},\boldsymbol{w})]\\=\sum_{k=1}^K \frac{|\Gamma^N_k |}{| \Gamma^N |}\max_{\mathbb{Q}_k\in \mathscr{W}_{\theta_k} }\mathbb{E}_{\mathbb{Q}_k}\bigg[ g(\boldsymbol{x},\boldsymbol{w})\bigg]\label{mid1}
\end{gather}

Based on this observation, we first derive equivalent reformulation of (\ref{BNWDRO}) for general loss function $g(\boldsymbol{x},\boldsymbol{w})$ in Theorem \ref{refor1}.
\begin{theorem}
    \label{refor1}
    Suppose that the loss function $g(\boldsymbol{x},\boldsymbol{w})$ is measurable. Then the BNWDRO problem (\ref{BNWDRO}) admits the following equivalent reformulation.
    \begin{equation}
        \begin{gathered}
                \label{measure}
        \min_{x\in \mathscr{X},\lambda_k\geq0,\alpha_{\widehat{\boldsymbol{w}}}}\left\{f(\boldsymbol{x})+\sum_{k=1}^K{\frac{|\Gamma^N_k |}{| \Gamma^N |}} \left[ \lambda_k \theta_k+\sum_{\widehat{\boldsymbol{w}}\in \Gamma_k^N}\frac{\alpha_{\widehat{\boldsymbol{w}}}}{|\Gamma^N_k|}  \right] \right\}\\
        \text{s.t.}\hspace{5pt} \alpha_{\widehat{\boldsymbol{w}}}\geq \max_{\boldsymbol{w}\in \varXi} \big\{g(\boldsymbol{x},\boldsymbol{w})-\lambda_k \| \widehat{\boldsymbol{w}}-\boldsymbol{w} \|\big\},\forall \widehat{\boldsymbol{w}}\in \Gamma^N
        \end{gathered}
    \end{equation}
\end{theorem}
\vspace{5pt}
\begin{proof}
    {According to Theorem 1 in \mbox{\cite{gao2023distributionally}}, each worst-case expectation in (\mbox{\ref{mid1}}) admits the following reformulation.}
    \begin{gather}
    {\max_{\mathbb{Q}_k\in \mathscr{W}_{\theta_k} }\mathbb{E}_{\mathbb{Q}_k}\bigg[ g(\boldsymbol{x},\boldsymbol{w})\bigg]}\nonumber\\
         {=\inf_{\lambda_k\geq0} \bigg\{ \lambda\theta+ \frac{1}{|\Gamma^N_k|}\sum_{i=1}^N \Big[  \max_{\boldsymbol{w}\in \varXi} \{ g(\boldsymbol{x},\boldsymbol{w})-\lambda_k \|  \widehat{\boldsymbol{w}}_i- \boldsymbol{w} \| \} \Big]           \bigg\}}\nonumber
    \end{gather}

   { By applying the above reformulation for each worst-case expectation term, (\mbox{\ref{measure}}) is derived immediately.}
\end{proof}

Since the reformulation (\ref{measure}) is still intractable due to the robust constraints, we pose restrictions to the support and loss functions to develop a tractable reformulation.

\begin{assumption}
\label{ass2}
    The feasible region $\mathscr{X}$ and uncertainty support $\varXi$ are a convex and closed, the loss function $f(\boldsymbol{x})$ is convex, and the loss function $g(\boldsymbol{x},\boldsymbol{w})=\max_{i\in [I]} h_i(\boldsymbol{x},\boldsymbol{w})$, where function $h_i(\boldsymbol{x},\boldsymbol{w})$ is proper, upper semicontinuous, convex in $\boldsymbol{x}$, and concave in $\boldsymbol{w}$ for each $i\in [I]$.
\end{assumption}

\begin{theorem}
\label{refor2}
    (Tractable reformulation) Suppose that Assumption \ref{ass2} holds. Then the BNWDRO problem (\ref{BNWDRO}) admits the following equivalent finite convex reformulation.
    \begin{equation}
    \label{trac}\nonumber
        \begin{gathered}
\min_{x,\lambda_k,\alpha_{\widehat{\boldsymbol{w}}},\boldsymbol{z}_{\widehat{\boldsymbol{w}},i},\boldsymbol{v}_{\widehat{\boldsymbol{w}},i}}\left\{f(\boldsymbol{x})+ {\small  \sum_{k=1}^K{\frac{|\Gamma^N_k |}{| \Gamma^N |}}}\left[ \lambda_k \theta_k+\sum_{\widehat{\boldsymbol{w}}\in \Gamma_k^N}\frac{\alpha_{\widehat{\boldsymbol{w}}}}{|\Gamma^N_k|}  \right] \right\}\\
            \text{s.t. } \alpha_{\widehat{\boldsymbol{w}}}\geq [-h_i]^\ast (\boldsymbol{x},\boldsymbol{z}_{\widehat{\boldsymbol{w}},i}-\boldsymbol{v}_{\widehat{\boldsymbol{w}},i})+\sigma_{\varXi}(\boldsymbol{v}_{\widehat{\boldsymbol{w}},i})-\langle \boldsymbol{z}_{\widehat{\boldsymbol{w}},i},\widehat{\boldsymbol{w}}\rangle,\\
            \forall \widehat{\boldsymbol{w}}\in \Gamma^N, \forall i\in [I]\\
            \| \boldsymbol{z}_{\widehat{\boldsymbol{w}},i} \|_{\ast}\leq \lambda_k,\forall k\in [K],\forall \widehat{\boldsymbol{w}}\in \Gamma_k^N, \forall i\in [I]\\
            \boldsymbol{x}\in \mathscr{X},\lambda_k\geq 0,\boldsymbol{z}_{\widehat{\boldsymbol{w}},i}\in \mathbb{R}^m,\boldsymbol{v}_{\widehat{\boldsymbol{w}},i}\in \mathbb{R}^m,\\
            \forall k\in [K],\forall \widehat{\boldsymbol{w}}\in \Gamma^N, \forall i\in [I]
        \end{gathered}
    \end{equation}
    where $[-h_i]^\ast(\boldsymbol{x},\boldsymbol{z})=\max_{\boldsymbol{w}\in \varXi}\{ \langle \boldsymbol{z},\boldsymbol{w}\rangle+h_i(\boldsymbol{x},\boldsymbol{w}) \}$ is the conjugate function of $-h_i$, $\| \cdot \|_{\ast}$ is the dual norm of $\| \cdot \|$, and $\sigma_{\varXi}$ is the support function of $\varXi$.
\end{theorem}
\begin{proof}
    By applying Theorem 4.2 in \cite{mohajerin2018data} to each worst-case expectation term with respect to the local Wasserstein ball in (\ref{mid1}), we can reach to the reformulation (\ref{trac}). 
\end{proof}

In many practical scenarios, the support is a convex polytope and the loss function $g(\boldsymbol{x},\boldsymbol{w})$ is convex piecewise linear. By leveraging Theorem \ref{refor2},  we present the reformulation for this special case in Corollary \ref{coro}.
\begin{corollary}
\label{coro}
    Suppose that Assumption \ref{ass2} holds with $h_i(\boldsymbol{x},\boldsymbol{w})$ affine in $\boldsymbol{w}$, \emph{i.e.}, $h_i(\boldsymbol{x},\boldsymbol{w})=\langle \boldsymbol{a}_i(\boldsymbol{x}),\boldsymbol{w}\rangle+b_i$, and that the support $\varXi=\{ \boldsymbol{C}\boldsymbol{w}\leq \boldsymbol{d} \}$ for some matrix $\boldsymbol{C}$ and vector $\boldsymbol{d}$. The BNWDRO problem (\ref{BNWDRO}) can be reformulated as follows.
        \begin{gather}
                    \min_{\boldsymbol{x},\lambda_k,\alpha_{\widehat{\boldsymbol{w}}},\boldsymbol{\psi}_{\widehat{\boldsymbol{w}},i}}\Bigg\{ f(\boldsymbol{x})+\sum_{k=1}^K \frac{|\Gamma^N_k |}{| \Gamma^N |}\bigg\{\lambda_k \theta_k+\sum_{\widehat{\boldsymbol{w}}\in \Gamma^N_k} \frac{\alpha_{\widehat{\boldsymbol{w}}}}{|\Gamma^N_k|} \bigg\}\Bigg\}\nonumber\\
        \text{s.t. }\alpha_{\widehat{\boldsymbol{w}}}\geq \langle \boldsymbol{a}_i(\boldsymbol{x}),\widehat{\boldsymbol{w}}\rangle+\langle \boldsymbol{\psi}_{\widehat{\boldsymbol{w}},i},\boldsymbol{d}-\boldsymbol{C}\widehat{\boldsymbol{w}}\rangle+b_i,\nonumber\\\forall \widehat{\boldsymbol{w}}\in\Gamma^N,\forall i\in [I]\nonumber\\
        \Big\|\boldsymbol{C}^T\boldsymbol{\psi}_{\widehat{\boldsymbol{w}},i}-\boldsymbol{a}_i(\boldsymbol{x})\Big\|_{\ast}\leq\lambda_k,\nonumber\\\forall k\in [K],\forall \widehat{\boldsymbol{w}}\in \Gamma_k^N, \forall i\in [I]\nonumber\\
        \boldsymbol{x}\in\mathscr{X},\lambda_k\geq 0,\boldsymbol{\psi}_{\widehat{\boldsymbol{w}},i}\geq 0, \forall k\in [K],\forall \widehat{\boldsymbol{w}}\in \Gamma^N, \forall i\in [I]\nonumber
        \end{gather}
\end{corollary}

\section{Case Study}

To validate the effectiveness of the BNWDRO framework, we conduct experiments on the UC problem under wind power generation uncertainty. {The UC model is adopted from \mbox{\cite{zhu2019wasserstein}} and experiments are conducted on the IEEE six-bus system and the IEEE 30-bus system. This problem aims to optimize the day ahead UC scheduling and the affine dispatch participation factors. For more information on this problem and how it is related to the general form (\mbox{\ref{BNWDRO}}), we refer the readers to the supplementary material.}


To test the performance of the proposed BNWDRO framework, we compare it with two commonly adopted DRO methods in the literature. The methods analyzed in this section are listed below.
\begin{itemize}
  \allowdisplaybreaks
  \item MDRO: DRO method using a moment ambiguity set \cite{zhou2019distributionally}, \emph{i.e.},
\begin{small}
  \begin{equation*}
    \mathscr{B}_\text{MDRO}=\left\{ \mathbb{Q} \in \mathscr{P}(\varXi)\left| \begin{gathered} \mathbb{E}_{\mathbb{Q}}(\boldsymbol{w})=\boldsymbol{\widehat{\mu}}
        \\\mathbb{E}_{\mathbb{Q}}(\boldsymbol{w}-\boldsymbol{\widehat{\mu}})^T(\boldsymbol{w}-\boldsymbol{\widehat{\mu}})\preceq \boldsymbol{\widehat{\Sigma}} \end{gathered}   \right. \right\}
  \end{equation*}  
\end{small}

\item WDRO: DRO method using the Wasserstein ambiguity set \cite{zhu2019wasserstein}.

\item BNWDRO: This work.
\end{itemize}

In determining the parameters of these DRO methods, we take $\widehat{\boldsymbol{\mu}}$ and $\widehat{\boldsymbol{\Sigma}}$ to be the mean value and covariance of the historical data in the MDRO method, respectively.
For the WDRO method and the proposed BNWDRO method, we determine their Wasserstein radius according to the following formula,  which is commonly adopted in power system applications \cite{zhu2019wasserstein}.
\begin{equation}
    \label{was-ra}
    \begin{gathered}
        \theta_k=C_k\sqrt{\frac{1}{ |\Gamma^N_k| }\text{log}\left( \frac{1}{1-\beta} \right)}\\
        C_k=2\inf_{z\geq0} \sqrt{\frac{1}{2z}\Bigg\{ 1+\text{ln}\Bigg( \frac{1}{|\Gamma^N_k|}{\footnotesize \sum_{\widehat{\boldsymbol{w}}\in \Gamma^N_k}}
        e^{z\|\widehat{\boldsymbol{w}} -\widehat{\boldsymbol{\mu}}_k \|_1^2} \Bigg) \Bigg\}}
    \end{gathered}
\end{equation}
where $\beta$ is the confidence level and $\widehat{\boldsymbol{\mu}}_k$ is the mean value of $\widehat{\boldsymbol{w}}\in \Gamma^N_k$. In all cases, the confidence level $\beta$ is set to 95\%.






In simulations, we apply each of these methods to 200 randomly generated data sets, which contain $N$ data independently and identically drawn from {the true wind power distribution $\mathcal{P}$. In computing the out-of-sample cost, we use the Montecarlo method with $1\times 10^6$ samples.}
All experiments are conducted on a laptop with Intel i5 CPU and 32G memory. All optimization methods are
established by YALMIP in MATLAB R2022a and solved by GUROBI 9.5.2.



{The simulation results of these DRO methods under different sample sizes on the IEEE 30-bus system are presented in Fig. \mbox{\ref{outofsample}}, and results on the IEEE six-bus system can be found in the supplementary material}. Fig. \ref{outofsample}(a) presents the certificated cost, \emph{i.e.}, the optimal value $\widehat{J}_N$ given by different methods. {It is noteworthy that the proposed BNWDRO method always gives a less conservative mean certificate cost compared with the other two methods.} 
Fig. \ref{outofsample}(b) depicts the out-of-sample performance of these methods. {It shows that the MDRO method fails to be asymptotically consistent.
For the two methods employing the Wasserstein metric, their out-of-sample costs consistently decrease as the data size $N$ grows and finally converge toward the optimal solution. However, when data is limited, the WDRO method proves excessively robust. On the contrary, the proposed BNWDRO is proved to be much less conservative.
More precisely, we conclude that compared with the MDRO and WDRO methods, the BNWDRO method
 demonstrates maximum improvements of 6,97\% and 7.34\% in out-of-sample mean cost, respectively. Fig. \mbox{\ref{outofsample}}(c) illustrates the reliability of the certificate cost, which shows that the MDRO and WDRO method is presented to be over-conservative. As for the proposed BNWDRO method, its reliability grows from 77\% as the sample size increases and converges to around 95\%, which is our prescribed reliability.}

\begin{figure}[htbp]
  \vspace{-7pt}
  \centering
    \includegraphics[width=0.8\hsize]{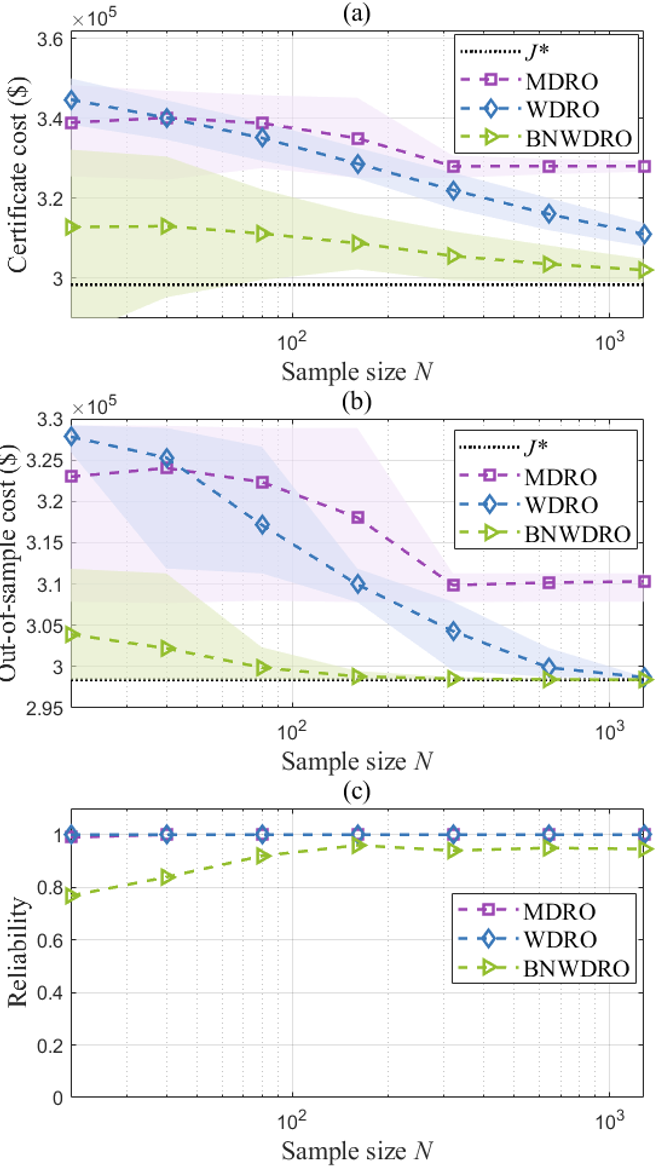}
    \caption{Simulation results under different sample sizes. (Results between the 10\% and 90\% 
quantiles are visualized in shade, the dashed lines stand for the mean values, and the dotted line represents 
the optimal cost derived from the SP counterpart under the true distribution.)}
    \label{outofsample} \vspace{-10pt}
    \end{figure}




\vspace{-6pt}
\section{Conclusion}
This work proposed a novel data-driven BNWDRO framework for decision-making under uncertainty. By unifying a Bayesian nonparametric method and the Wasserstein metric, the proposed BNWDRO framework can fully explore and characterize the global-local features of the data, thereby providing less conservative decisions for the decision-maker. By investigating the connection with the Wasserstein DRO, we proved that the proposed BNWDRO framework enjoys theoretical guarantees.
Subsequently, we derived an equivalent reformulation to solve the resulting DRO problem under the proposed framework and verified the effectiveness of the BNWDRO framework on a UC problem.

\vspace{-6pt}
\bibliographystyle{IEEEtran}
\bibliography{ref.bib}

\end{document}


\clearpage\setlength{\parskip}{.3em}
{\centering \bf \large Supplementary Material of "Data-Driven Bayesian Nonparametric Wasserstein
Distributionally Robust Optimization" \par}
\hfill\\
{\centering Chao Ning and Xutao Ma \par}
\hfill

\section{Detailed Unit Commitment Model}

We adopt The UC model proposed in \cite{zhu2019wasserstein} and we present the objective function of this model below.
    \allowdisplaybreaks
    \begin{gather}
    \label{r1obj}
        \begin{gathered}
        \min_{(s_{it},u_{it},w_{it},\overline{r}_{it},\underline{r}_{it},P_{it}^g,\alpha_{it}:\ t\in \mathcal{T},i\in \mathcal{G})\in \mathscr{X}} \bigg\{  \sum_{t\in \mathcal{T}}\sum_{i\in \mathcal{G}} \big(SU_i u_{it}+SD_i w_{it}+CU_i \overline{r}_{it}+CD_i\underline{r}_{it}\\+F_{it}(P^{g}_{it})\big)+\sup_{\mathbb{Q}\in \mathscr{B}} \mathbb{E}_{w\sim \mathbb{Q}}\Big[\sum_{t\in \mathcal{T}}\sum_{i\in \mathcal{G}} \alpha_{it} d_i | w | \Big]    \bigg\}
        \end{gathered}
    \end{gather}
    where $F_{it}(\cdot)$ is a piece-wise linear function and the feasible region $\mathscr{X}$ is given by
    \begin{subequations}
        \begin{gather}
        \mathscr{X}=\{ (s_{it},u_{it},w_{it},\overline{r}_{it},\underline{r}_{it},P_{it}^g,\alpha_{it}:\ t\in \mathcal{T},i\in \mathcal{G})|\nonumber\\
        TU_t(s_{it}-s_{i(t-1)})\leq \sum_{k=0}^{TU_i-1}s_{it+k},\ \forall t\in \mathcal{T},\forall i\in \mathcal{G}\label{r1con1} \\
        TD_i(s_{it}-s_{i(t-1)})\leq \sum_{k=0}^{TD_i-1}(1-s_{it+k}),\ \forall t\in \mathcal{T},\forall i\in \mathcal{G}\\
        s_{it}-s_{i(t-1)}=u_{it}-w_{it},\ \forall t\in \mathcal{T},\forall i\in \mathcal{G}\label{r1con2}\\
        \sum_{i\in \mathcal{G}}P^g_{it}+\sum_{j\in \mathcal{W}}P_{jt}^w=P_t^l,\ \forall t\in \mathcal{T}\label{r1con3}\\
        s_{it} \underline{P}_i^g+\underline{r}_{it}\leq P_{it}^g\leq s_{it} \overline{P}^g_i-\overline{r}_{it},\ \forall t\in \mathcal{T},\forall i\in \mathcal{G}\label{r1con4}\\
        \begin{gathered}
            (P^g_{it}+\overline{r}_{it})-(P^g_{i(t-1)}-\underline{r}_{i(t-1)})\leq (2-s_{i(t-1)}-s_{it})\overline{RU}_i+(1+s_{i(t-1)}-s_{it})RU_i\\,\ \forall t\in \mathcal{T},\forall i\in \mathcal{G}
        \end{gathered}\label{r1con5}
        \\
        \begin{gathered}
           (P^g_{i(t-1)}+\overline{r}_{i(t-1)})-(P^g_{it}-\underline{r}_{it})\leq (2-s_{i(t-1)}-s_{it})\overline{RD}_i+(1-s_{i(t-1)}+s_{it})RD_i \\,\ \forall t\in \mathcal{T},\forall i\in \mathcal{G}
        \end{gathered}\label{r1con6}
        \\
        \sum_{i\in\mathcal{G}}\alpha_{it}=1,\ \forall t\in \mathcal{T}\label{r1con7}\\
        0\leq \alpha_{it}\leq 1,\ \forall t\in \mathcal{T},\forall i\in \mathcal{G}\\
        -\underline{r}_{it}\leq -\alpha_{it}\overline{w},-\alpha_{it}\underline{w}\leq\overline{r}_{it},\ \forall t\in \mathcal{T},\forall i\in \mathcal{G} \}\label{r1con8}
    \end{gather}
    \end{subequations}
    where (\ref{r1con1})-(\ref{r1con2}) is the physical constraints on the unit on/off state, (\ref{r1con3}) is the power balance constraint, (\ref{r1con4}) restricts the minimum and maximum power output of unit $i$ at time $t$, (\ref{r1con5})-(\ref{r1con6}) set limitations on the unit ramping up and ramping down power and (\ref{r1con7})-(\ref{r1con8}) correspond to the constraints of affine dispatch policy.

    The uncertainty of this UC model is the total wind power forecasting error $w$, whose distributional ambiguity is captured in the ambiguity set $\mathscr{B}$. The objective function (\ref{r1obj}) aims to minimize the sum of start-up cost $SU_i u_{it}$, shut-down cost $SD_i w_{it}$, upward reserve cost $CU_i \overline{r}_{it}$, downward reserve cost $CD_i\underline{r}_{it}$, generation cost $F_{it}(P^{g}_{it})$ and the worst-case expected adjusting cost $\sup_{\mathbb{Q}\in \mathscr{B}} \mathbb{E}_{w\sim \mathbb{Q}}\Big[\sum_{t\in \mathcal{T}}\sum_{i\in \mathcal{G}} \alpha_{it} d_i | w | \Big] $. In relation to the general form 
    \begin{equation}
        \min_{x\in \mathscr{X}} \bigg\{ f(\boldsymbol{x})+\max_{\mathbb{Q}\in \mathscr{B}(\Gamma^N)}\mathbb{E}_{\mathbb{Q}}[g(\boldsymbol{x},\boldsymbol{w})]\bigg\} ,
    \end{equation}
    we have
    \begin{gather}
    \boldsymbol{x}=(s_{it},u_{it},w_{it},\overline{r}_{it},\underline{r}_{it},P_{it}^g,\alpha_{it}:\ t\in \mathcal{T},i\in \mathcal{G})\\
        f(\boldsymbol{x})=\sum_{t\in \mathcal{T}}\sum_{i\in \mathcal{G}} \big(SU_i u_{it}+SD_i w_{it}+CU_i \overline{r}_{it}+CD_i\underline{r}_{it}+F_{it}(P^{g}_{it})\big)\\
        g(\boldsymbol{x},\boldsymbol{w})=\max \Bigg\{\sum_{t\in \mathcal{T}}\sum_{i\in \mathcal{G}} \alpha_{it} d_i  w,-\sum_{t\in \mathcal{T}}\sum_{i\in \mathcal{G}} \alpha_{it} d_i  w\Bigg\}
    \end{gather}
    
    Note that $g(\boldsymbol{x},\boldsymbol{w})$ in this UC model satisfies the condition in Corollary 1, so it can be reformulated into a mixed-integer linear program.

    \begin{figure}[htbp]
  \vspace{0pt}
  \centering
    \includegraphics[width=0.55\hsize]{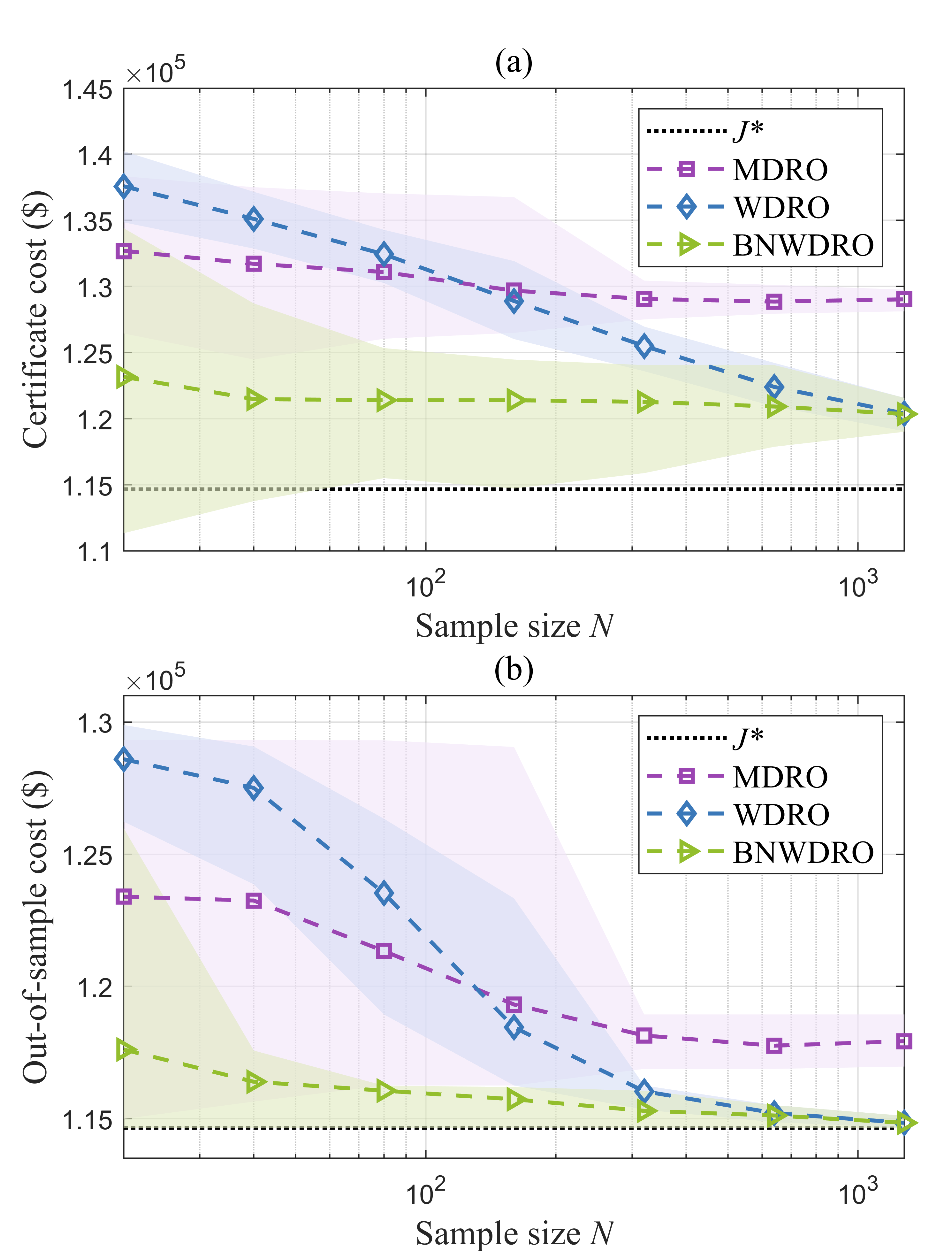}
    \caption{Out-of-sample performance and certificate under different sample sizes. (Results between the 10\% and 90\% 
quantiles are visualized in shade, the dashed lines stand for the mean values, and the dotted line represents 
the optimal cost derived from the SP counterpart under the true distribution.)}
    \label{outofsample}
    \end{figure}

\section{Results on the IEEE six-bus system}

The simulation results of these DRO methods under different sample sizes are presented in Fig. \ref{outofsample}. Fig. \ref{outofsample}(a) presents the certificated cost, \emph{i.e.}, the optimal value $\widehat{J}_N$ given by different methods. It is noteworthy that by incorporating the global-local information, the proposed BNWDRO method always gives less conservative certificate cost compared with the other two methods. 
Fig. \ref{outofsample}(b) depicts the out-of-sample performance of these methods. It shows that the MDRO method fails to be asymptotically consistent and converges to a suboptimal value as the sample size grows.
For the two methods employing the Wasserstein metric, their out-of-sample costs consistently decrease as the 
data size $N$ increases. Notably, when a substantial amount of data becomes available, these two methods all converge toward the optimal solution. However, in scenarios where data is limited, the WDRO method proves excessively robust. Fortunately, the unification with the DPMM can remarkably mitigate
conservatism, as observed in the BNWDRO method. 
This observation underscores the effectiveness of deciphering and incorporating the
{global-local features of data}. 
More precisely, we conclude that compared with MDRO and WDRO methods, the BNWDRO method
 demonstrates maximum improvements of 4.70\%, and
8.55\% in out-of-sample mean cost, respectively.

\bibliographystyle{unsrt}
\bibliography{ref.bib}